\documentclass[12pt,a4paper]{article}
\usepackage{latexsym,amssymb,amsmath,amsthm,amsfonts,euscript}
\usepackage[cp1251]{inputenc}

\newcommand{\ignore}[1]{}

\newcommand{\B}{\mathcal{B}}

\newcommand{\ZZ}{\mathbb{Z}}

\newcommand{\ind}{{\rm ind}}

\begin{document}
\title{An explicit finite $B_k$-sequence}
\date{}
\author{Igor S. Sergeev\footnote{e-mail: isserg@gmail.com}}

\maketitle

\begin{abstract}
For any $n$ and $k$, we provide an explicit (that is, computable
in polynomial time) example of integer $\B_k$-sequence of size $n$
consisting of elements bounded by $n^{k+o(k)}$.
\end{abstract}

\begin{flushright}
{\it dedicated to the memory of Vladimir Evgen'evich Alekseev} (1943--2020)
\end{flushright}

{\bf Introduction.} Recall that a set $B$ in some commutative
group is a {\it $\B_k$-sequence} if all $k$-element sums in $B$
are different, that is, the equality
\[  a_1 + \ldots + a_k = b_1 + \ldots + b_k, \qquad a_i, b_j \in B, \]
holds iff the multisets of summands coincide: $\{a_1, \ldots,
a_k\} = \{b_1, \ldots, b_k\}$.

$\B_2$-sequences are also known as {\it Sidon sequences}. Very
often, the notion of Sidon sequence stands as a synonym for
$\B_k$-sequence in general.

Easy to check, if $\ZZ_N$ contains a size-$n$ $\B_k$-sequence,
then $N \ge \binom{n+k-1}{k}$. We want to consider only
satisfactorily dense size-$n$ $\B_k$-sequences, say, for $N =
(n+k)^{O(k)}$, avoiding trivial examples like $\{k, k^2, \ldots
k^n\}$ with exponentially large elements. Also, we interest in
{\it explicit} constructions, that is, those that can be computed
in polynomial time with respect to the binary size\footnote{That
is, the length of the binary code representing the elements of the
set.}.

{\bf History.} The most famous explicit examples of the optimal
density integer Sidon sequences are: a size-$(q+1)$ set in
$\ZZ_{q^2+q+1}$ due to J.~Singer~\cite{sin38}, a size-$q$ set in
$\ZZ_{q^2-1}$ due to R.~C.~Bose~\cite{bos42}, and a size-$(p-1)$
set in $\ZZ_{p^2-p}$ due to V.~E.~Alekseev~\cite{ale81}. Here $p$
and $q$ stay for any prime number and prime power, respectively.
The latter set is attributed to I.~Ruzsa~\cite{ruz93} almost
everywhere.

The classical example of a nearly dense-optimal $\B_k$-sequence
was proposed by Bose and S.~Chowla in~\cite{bc62}. Let us recall
this construction that generalizes~\cite{bos42}. Let $GF(q) =
\{\alpha_1,\ldots,\alpha_q\}$, and $x$ be a primitive element in
$GF(q^k)$. It can be easily verified that the set
\[ D[q,k] = \{ d_i \mid x^{d_i} = x + \alpha_i,\; 1\le d_i < q^k \}\]
is a size-$q$ $\B_k$-sequence in $\ZZ_{q^k-1}$.

There are known also a number of similar constructions including another $\B_k$-sequence
from~\cite{bc62} generalizing~\cite{sin38}. 
H.~Derksen~\cite{der04} proposed even more general constructions
considering quotient polynomial rings $GF(q)[x]/(P(x))$ instead of
pure fields in the examples from~\cite{bc62}. C.~A.~G\'omez Ruiz
and C.~A.~Trujillo Solarte~\cite{gt11} extended an
example~\cite{ale81} to $\B_k$-sequences in $\ZZ_{p^k-p}$.

{\bf Discussion.} All these examples of $\B_k$-sequences may be
considered explicit only for constant or extremely slowly growing
$k$'s with respect to $n$, since they imply computation of
discrete logarithms in groups of generally non-smooth order.
Indeed, probabilistic or greedy constructions that we haven't
mentioned are even less explicit. It looks like we lack easily
computable and dense enough examples of $\B_k$-sequences that
could be useful in some specific situations, e.g. for proving
explicit lower bounds in computational complexity~\cite{is23}.
Thus, we intend to close this gap.

We follow the general idea of previous constructions: computing an
additive numeric $\B_k$-sequence as an image of some simple
multiplicative $\B_k$-sequence from an appropriate group. All we
need to make computations easy is to choose a basic multiplicative
group of smooth order. Note that in doing this, we will partially
sacrifice the density.

{\bf Construction.} Further, $p_1, p_2, \ldots$ denote odd prime
numbers written in growing order. Let $r = 1 + \lceil k\log p_n
\rceil$. The set of odd numbers-residues from 1 to $2^r-1$
constitutes the multiplicative group $\ZZ^*_{2^r}$ of the ring
$\ZZ_{2^r}$. For $r\ge3$, this group is a direct product of cyclic
groups of orders 2 and $2^{r-2}$, namely, $\ZZ^*_{2^r} \cong
\langle -1 \rangle_2 \langle 5 \rangle_{2^{r-2}}$ with $-1$ and
$5$ being generating elements. Therefore, any odd number $x$ has a
unique representation $x \equiv (-1)^j\cdot 5^h\; (\bmod\; 2^r)$,
where $0 \le j \le 1$ and $0 \le h < 2^{r-2}$. For details, see
e.g.~\cite{vin52}.

Consider the number set
\[ H[n,k] = \{ h_i \mid p_i \equiv \pm 5^{h_i}\; (\bmod\; 2^r),\, 0 \le h_i < 2^{r-2},\,i=1,\ldots,n\}.\]
Let us check that the given set is a $\B_k$-sequence in
$\ZZ_{2^{r-2}}$. By the choice of $r$, for different tuples of
indices $1 \le i_1 \le \ldots \le i_k \le n$, all numbers $\pm
p_{i_1} \cdot \ldots \cdot p_{i_k}$ are different and do not
exceed $2^{r-1}-1$ by absolute value. Hence, all residues
$5^{h_{i_1} +\, \ldots +\, h_{i_k}\;} (\bmod\; 2^r)$ are
different, and all sums $h_{i_1} +\, \ldots +\, h_{i_k\;} (\bmod\;
2^{r-2})$ are different as well.

The set $H[n,k]$ is not as dense as $D[q,k]$ or similar
constructions. Still, its density is satisfactorily in asymptotic
sense: $ 2^{r-2} < p_n^k < (2n\log (n+2))^k $ due to the
well-known facts about distribution of prime numbers, see
e.g.~\cite{rs62}.

We are left to confirm explicitness: that the set $H[n,k]$
requires ${(n+k)^{O(1)}}$ time to be constructed. First, we need
to obtain the list of prime numbers. Second, we have to compute
discrete logarithms\footnote{Here, we don't resort to the commonly
used notation $\ind_g\, x$.} $\log_5 (\pm p_i)$ in $\ZZ_{2^r}$.
For the first part, we may use Eratosthenes sieve or any other
known algorithm running in time $n^{O(1)}$. Discrete logarithm in
the cyclic group of order $2^{r-2}$ may be computed trivially by
$O(r^2)$ elementary arithmetic operations mostly consisting of
squarings. Indeed, we may determine binary digits of the number $a
= [a_{r-3},\, \ldots, a_0]_2 = \log_5 x\; (\bmod\; 2^r)$ 
sequentially as
\begin{multline*}
 a_0 = \log_{5^{2^{r-3}}} x^{2^{r-3}}, \quad a_1 = \log_{5^{2^{r-3}}}
 (5^{-a_0}x)^{2^{r-4}}, \quad  \ldots, \\
  a_{r-3} = \log_{5^{2^{r-3}}} \big(5^{-2^{r-4}a_{r-4} - \ldots - 2a_1 -
  a_0}x\big).
\end{multline*}
Inner logarithms are performed in an order-2 subgroup with
generating element $5^{2^{r-3}} \equiv 2^{r-1}+1\; (\bmod\; 2^r)$
simply by comparing with 1 and $2^{r-1}+1$. If both comparisons
fail, then $x \notin \langle 5 \rangle_{2^{r-2}}$.

{\bf Notes.} In the above example, we intentionally used as smooth
order for the basic multiplicative group as possible. Instead, we
can work in any ring $\ZZ_{p^r}$ with an odd prime $p$. The
multiplicative group $\ZZ^*_{p^r}$ has order $(p-1)p^{r-1}$ and it
is cyclic. The case $p=3$ is especially attractive, since there we
have 2 as a generating element for the multiplicative group. With
more care, we can consider residue rings of some other smooth
orders.

The choice of prime numbers for a ``factor base'' is also changeable.
Say, we can relax the condition of being prime to the
condition of being pairwise prime. Though, this relaxation alone
doesn't allow to substantially increase the density of the set.

Essentially, the present text in an excerpt from~\cite{is23}.

\end{document}